\documentclass[11pt,oneside,a4paper]{article}
\usepackage{syntonly}
\usepackage{geometry}
\usepackage{mathrsfs}
\usepackage{makeidx}
\makeindex
\usepackage{amsmath}
\usepackage{amssymb}
\usepackage{amsmath}
\usepackage{amsfonts}
\usepackage{amssymb}
\usepackage{esvect}
\usepackage{algorithm}
\usepackage{bm}
\usepackage{tikz}
\usepackage{cleveref}
\usepackage{graphicx}
\usepackage{epstopdf}

\newtheorem{proposition}{Proposition}

\newcommand{\qed}{\nobreak \ifvmode \relax \else
      \ifdim\lastskip<1.5em \hskip-\lastskip
      \hskip1.5em plus0em minus0.5em \fi \nobreak
      \vrule height0.30em width0.4em depth0.25em\fi}

 \author{Safari Mukeru and Mmboniseni P. Mulaudzi\\
\footnotesize{\em Department of Decision Sciences}\\
\footnotesize{University of South Africa, P. O. Box 392, Pretoria, 0003. South Africa}\\
\footnotesize{e-mails: mukers@unisa.ac.za, mulaump@unisa.ac.za}}

\title{{On the inverse problem of fractional Brownian motion and the inverse of infinite Toeplitz matrices}}

\date{}

\begin{document}

\maketitle

\pagenumbering{arabic}

\begin{abstract}

The inverse problem of fractional Brownian motion and other Gaussian processes with stationary increments involves inverting an infinite hermitian positively definite Toeplitz matrix (a matrix that has equal elements along its diagonals).
The problem of inverting Toeplitz matrices is  interesting on its own and has various applications in physics, signal processing, statistics, etc. A large body of literature has emerged to study this question since the seminal work of Szeg\"o on Toeplitz forms in 1920's. 
In this paper we obtain, for the first time, an explicit general formula for the inverse of infinite hermitian positive definite Toeplitz matrices. Our formula is explicitly given in terms of the Szeg\"o function associated to the spectral density of the matrix. These results are applied to the  fractional Brownian motion and to $m$-diagonal  Toeplitz matrices and we provide explicit examples. 

\end{abstract}

\section{Introduction}
The statistical inverse solution to the inverse problem involves quite often the inverse of a large or infinite matrix. Generally the measurement vector is modeled as a realisation of a stochastic process with unknown parameters and the problem is to recover these parameters from the measurements. In this paper we are mainly concerned with the fractional Brownian motion which is a Gaussian process $\{X(t): t\geq 0\}$ defined by the covariance function 
$$\langle X(t), X(s)\rangle =  \frac{a}{2} (t^{2H} + s^{2H} - |t-s|^{2H}),\,\,\,t, s \geq 0$$
where $H \in (0, 1)$ is a parameter called the Hurst index and $a >0$ is a scale parameter. The fractional Brownian motion has stationary increments but these increments are not independent. It is an important model that has been  successfully used to model  real data that appear in different contexts such as anomalous diffusion, single-file dynamics, volatility in financial markets, hydrology, geophysics, 
 queueing networks such as internet traffic, etc. (\cite{Cintoli_et_al}, \cite{Delorme_2015}, \cite{Rostek}, \cite{Vehel}).
 In what follows we shall assume that $a = 1$ and in this case $X(t)$ is called the standard fractional Brownian motion. 
The inverse problem of fractional Brownian motion deals with the recovery of the Hurst parameter $H$  from a given vector of measurements. For this it is necessary to consider the increments $\Delta_n = X(n+1) - X(n)$, $n \in \mathbb{N}$, of the fractional Brownian motion $X(t)$ and the corresponding covariance matrix $G = (g_{n,m})$ given by 
 \begin{eqnarray} \label{sddsd3eres}
g_{n, m} = \langle \Delta_n, \Delta_m \rangle =  {\scriptstyle\frac{1}{2}}|n-m+1|^{2H} + {\scriptstyle\frac{1}{2}}|n-m-1|^{2H} -|n-m|^{2H},\,\,\,n,m = 1,2,3\ldots
\end{eqnarray}
The process $\{\Delta_n: n \in \mathbb{N}\}$ is called the (discrete) fraction Gaussian noise. The statistical inverse method for fractional Brownian motion requires the computation of the inverse of the matrix $G$. In fact assume that $X_n = X(t_1), X(t_2), \ldots, X(t_n)$ is the vector of the observed values of the process at times $0\leq t_1 < t_2 < \ldots t_n \leq 1$. Take $t_i = i/n$ for all $i$. Set $W_n = (W_1, W_2, \ldots, W_n)$ where
 $W_i = X(t_{i}) - X(t_{i-1})$, $i = 1,2,\ldots, n$ with $X(t_0) = 0$.   The unknown parameter $H$ is considered as a random variable.  The conditional distribution of the random vector $W_n$ under the condition that the parameter $H$ is equal to $\hat H$ is: 
    $$D(W_n|\hat H) = \frac{D_{post}(\hat H, W_n)}{D_{prior}(\hat H)}.$$
The joint distribution $D_{post}(\hat H, W_n)$ of $\hat H$ and $W_n$ is given by:
    
$$D_{post}(\hat H, W_n) = \frac{1}{(2\pi)^{n/2} \left(\det(M_n)\right)^{1/2}} \exp\left(-\scriptstyle\frac{1}{2} W_n^T (M_n)^{-1} W_n\right)$$ where $M_n$ is the covariance matrix of $W_n$ (increments of the fractional Brownian motion with parameter $\hat H$). Clearly, 
        $$M_n = \frac{1}{n^{2\hat H}} G_{n\times n}$$ where $G_{n\times n}$ is the covariance matrix of the fractional Gaussian noise $(\Delta_1, \Delta_2, \ldots, \Delta_n)$. Of course it is necessary to take $n$ very large for a better approximation and ultimately $n \to \infty$  for the exact value.  However the value of  $(G_{n \times n})^{-1}$ for large $n$ and  the value of  $G^{-1}$ are not now known. 
 Up to so far in the literature, because $G^{-1}$ is unknown, it has to be approximated by some other infinite matrices.  One such approximation is to consider the inverse of the matrix $\tilde G = (\tilde g_{n,m})$ given by 
            $$\tilde g_{n,m} = g_{n,m} \mbox{ for } n - m \in \{ 0, 1, -1\} \mbox{ and } 0 \mbox{ otherwise.}$$
 (See for example \cite{Dambrogi-ola}.) A more widely used approximation in the literature is the asymptotic approximation introduced by Whittle  \cite{Whittle} in 1953. It consists of the matrix $\Gamma$ given by:
   $$ \Gamma_{k,j} = \int_0^1 \frac{e^{-2 \pi i (k-j)t}}{\varphi_H(t)}\,dt,\,\,k,j = 1,2,3, \ldots$$ 
   where $\varphi_H(t)$ is the spectral density function of the fractional Gaussian noise of index $H$ (introduced in section \ref{spectral}). (We refer to Beran \cite{Beran} for more details.)  This matrix $\Gamma$ corresponds to the inverse $G^{-1}$ only asymptotically in the sense that
 $$\lim_{k\to \infty} \Gamma_{k,j} = \lim_{k \to \infty} \left(G^{-1}\right)_{k,j} \,\,\, \mbox{ for each } j \mbox{ fixed}.$$
One of the objectives of this paper is to complete these earlier efforts of Whittle and others by providing for the first time an analytic expression that gives the exact value of the matrix $G^{-1}$. 

Our method is robust in the sense that it applies not only to the covariance matrix $G$ of the fractional Gaussian noise but generally to other infinite hermitian positive definite Toeplitz matrices. In what follows $G$ will be an infinite hermitian Toeplitz matrix on the set of complex numbers. 

The problem of inverting infinite Toeplitz matrices has important applications in many  contexts and it is interesting on its own (see for example \cite{Caflisch},  \cite{daFonseca}, \cite{Pedro}, \cite{Wang}).   This work is a novel contribution to this problem.

 \section{Spectral density function} \label{spectral}
Let $G = (g_{k,j})_{k,j=1,2,\ldots}$ be an infinite, Hermitian and positive definite Toeplitz matrix. We assume that $G$ is such that 
$$g_{k,j} = \gamma(j-k),\,\,k,j = 1,2,3,\ldots$$ where $\gamma$ is a positive definite function defined on the set of integers $\mathbb{Z}$ and taking values in the set of complex numbers $\mathbb{C}$ and satisfies
$$\gamma(0) = 1 \mbox{ and } \gamma(-k) = \overline{\gamma(k)}, \mbox{ for all } k \in \mathbb{Z}.$$
This implies (by the classical Bochner theorem) that there exists a probability measure $\mu$ on the unit interval $[0, 1]$  such that
$$\gamma(k) = \int_0^1 e^{2  \pi k i t} d\mu(t)$$ for all $k \in \mathbb{Z}$. We shall assume that $\mu$ admits a probability density function $\varphi \geq 0$ that is integrable (more precisely $\varphi^p$ is integrable  for some number $p> 1$) so that its Fourier series 
\begin{eqnarray} \label{sdsde34re}
\varphi(t) \sim \sum_{k\in \mathbb{Z}} \gamma(k) e^{-2 \pi k i t}
\end{eqnarray}
 converges for almost all $t \in [0, 1]$ (with respect to the Lebesgue measure on $[0, 1]$).
 The function $\varphi$ is the spectral density function of the matrix $G$. We shall assume that $G$ is invertible in the sense that there exists another positive definite matrix $G^{-1}$ such that 
$G G^{-1} = G^{-1} G = I$ where $I$ is the infinite identity matrix. Our goal is to compute the inverse matrix $G^{-1}$.  

If $G$ is the covariance matrix of  fractional Gaussian noise given by relation (\ref{sddsd3eres}), then clearly $G$ is the Toeplitz matrix corresponding to the function 
$$\gamma(k) =  {\scriptstyle\frac{1}{2}}|k+1|^{2H} + {\scriptstyle\frac{1}{2}}|k-1|^{2H} -|k|^{2H}.$$ It is an important result obtained by Sinai \cite{Sinai} that $G$ admits a spectral density function $\varphi_H$ given by 
 \begin{eqnarray*} \label{dsdfesw23}
\varphi_H(t) = C(H)|e^{2 \pi i t}-1|^2\left(\sum_{n=-\infty}^\infty \frac{1}{|t+n|^{2H+1}}\right),\, 0 < t < 1. 
\end{eqnarray*}
Here $C(H)$ is a normalising constant given by $$C(H) = -\zeta(-2H)/(2 \zeta(1+2H))$$ where $\zeta(.)$ is the Riemann zeta function.  Clearly, 
\begin{eqnarray*} \label{dsdfesw2312s}
\varphi_H(t) &=& 4 C(H) \left(\sin^2 \pi t\right) \sum_{n=0}^\infty \left(\frac{1}{(n+t)^{2H+1}} + \frac{1}{(n+ 1-t)^{2H+1}}\right)\\
            & = &4 C(H) \left(\sin^2 \pi t\right) (\zeta(2H+1, t) + \zeta(2H+1, 1-t))
\end{eqnarray*}
 where $\zeta(.,.)$ is the classical Hurwitz zeta function. 
Moreover it is well-known that 
$$\varphi_H(t)  = O(t^{1-2H} (1-t)^{1-2H})$$  for  $ t $   near 0 or near  1. (See \cite{Mukeru_Pisieral} for details.) 
This implies that for $0 < H \leq 1/2$, the function  $\varphi_H$ is continuous on the interval $[0, 1]$ and for $1/2 < H < 1$, it is not continuous but $(\varphi_H)^p$ is integrable on the interval $[0, 1]$ for all $1 < p <1/(2H-1)$. The Fourier series of $\varphi_H$ (given by relation (\ref{sdsde34re})) yields the following representation: 
\begin{eqnarray*}
\varphi_H(t)= 
 (\cos(2\pi t)-1)\left(\mbox{Li}(-2 H, e^{-2\pi i t}) + \mbox{Li}(-2 H, e^{2\pi i t})\right),\,\,0 < t < 1, 
\end{eqnarray*}
where $\mbox{Li}$ is the analytic continuation of the classical polylogarithm function 
 $\mbox{Li}(s,z) = \sum_{k=1}^\infty  k^{-s} z^k $ in the complex plane (except at the single point $z = 1$). 
This representation is amenable to calculations (for example the function $\mbox{Li}$ is implemented in Wolfram Mathematica).

\section{Szeg\"o function of a Toeplitz matrix}
Assume that the matrix $G$ admits a spectral density function $\varphi$. Under the condition that the spectral density function $\varphi$ satisfies the Szeg\"o condition
\begin{eqnarray} \label{Szego_cond}
\int_0^1 \log(\varphi(t)) dt > -\infty, 
\end{eqnarray}
 one can associate to the matrix $G$ the function $S: D \to \mathbb{C}$ (where $D = \{z \in \mathbb{C}: |z| < 1\}$ is the unit disc) defined by
        \begin{eqnarray} \label{dsdferef3}
S(z) = \exp\left({\scriptstyle\frac{1}{2}}\int_0^1 \left(\frac{e^{2\pi i t} + z}{e^{2\pi i t} - z}\right) \log(\varphi(t)) dt\right),\,\, \mbox{ for all } |z| < 1. 
\end{eqnarray}
The function $S(z)$ is known as the Szeg\"o function associated to the matrix $G$ (or equivalently associated to $\varphi$). (We refer the reader to Grenander and Szeg\"o \cite{Szego} and Simon \cite{Simon1} for more details.)

We shall mainly use the inverse $1/S(z)$ and we set 
  \begin{eqnarray} \label{dsdferef3qw23}
\psi(z) = S^{-1}(z) =  \exp\left({-\scriptstyle\frac{1}{2}}\int_0^1 \left(\frac{e^{2\pi i t} + z}{e^{2\pi i t} - z}\right) \log(\varphi(t)) dt\right),\,\, \mbox{ for all } |z| < 1. 
\end{eqnarray}
It is also known that if the Szeg\"o function $S(z)$ is written in the form,
  \begin{eqnarray} \label{szego1}
S(z) = c_0 + c_1 z + c_2 z^2 + \ldots,\,\,|z| < 1
\end{eqnarray}
where $c_0, c_1, c_2,\ldots$ are complex numbers, then these coefficients  satisfy the system: 
\begin{eqnarray} \label{szego2}
\int_0^1 e^{-2\pi i k t} \varphi(t) dt = \overline{c_0} c_k + \overline{c_1} c_{k+1} + \overline{c_2} c_{k+2} + \ldots,\,\,k = 0,1,2,\ldots
\end{eqnarray}
 together with the condition that $c_0$ is a positive real number. 
 
Let 
 $P_1(z), P_2(z), P_3(z),\ldots$ be a sequence of polynomials of a complex variable $z$ which are orthonormal on the unit circle $z = e^{2\pi i t}$, $t\in [0, 1]$ with respect to the weight function $\varphi(t)$. That is $$\int_0^1 P_n(e^{2\pi i t}) \overline{P_m(e^{2\pi i t})} \varphi(t) dt = \delta_n^m,\,\,\, n,m = 1,2,3,\ldots$$
Assume moreover that  $\mbox{degree}(P_n(z)) = n-1$ for all $n = 1,2,\ldots$ and further that the leading coefficient of $z^{n-1}$ is a positive real number. These conditions uniquely determine the sequence $(P_n(z))$. 
Szeg\"o  proved that the function $\psi$ is such that for all $z, w \in D$, 
             $$\sum_{n= 1}^\infty P_n(z) \overline{P_n(w)} =\frac{ \psi(z) \overline{\psi(w)}}{1 - z \overline{w}}.$$
(See the book by Grenander and Szeg\"o \cite[pp. 37-51]{Szego}.)
          
For the case of the fractional Gaussian noise, since as mentioned before, the spectral density function 
 $\varphi_H(t) $ is continuous in the open interval $(0, 1)$ and satisfies 
  $$\varphi_H(t) = O(t^{1-2H} (1-t)^{1-2H})$$ at the boundary of the interval, that is for $t$ near 0 and $t$ near 1, then it is clear that  $\varphi_H$ satisfies the Szeg\"o condition (\ref{Szego_cond})  and hence the fractional Gaussian noise admits a Szeg\"o function. We shall denote its inverse by $\psi_H$. That is, 
  \begin{eqnarray} \label{dsdferef3232q}
\psi_H(z) = \exp\left({-\scriptstyle\frac{1}{2}}\int_0^1 \left(\frac{e^{2\pi i t} + z}{e^{2\pi i t} - z}\right) \log(\varphi_H(t)) dt\right),\,\, \mbox{ for all } |z| < 1.  
\end{eqnarray}

  \section{Szeg\"o function and the inverse of an infinite Toeplitz matrix }
Let $Q_1(z), Q_2(z), Q_3(z), \ldots$ be arbitrary sequence of orthonormal polynomials on the unit circle with respect to the weight function $\varphi$ such that $\mbox{deg}(Q_n) = n-1$. Here it is not necessary to impose  any extra condition on the coefficient of $z^{n-1}$ in $Q_n(z)$ as it is the case for the polynomials $P_n(z)$. 
We have that any such general sequence of orthogonal polynomials satisfies:
  $$\sum_{n= 1}^\infty Q_n(z) \overline{Q_n(w)} =\frac{ \psi(z) \overline{\psi(w)}}{1 - z \overline{w}}.$$
  This is a consequence of an important result about the structure of the space $\mathscr{H}$ of holomorphic functions $f: D \to \mathbb{C}$ of the form:
    $$f(z) = \sum_{n=1}^\infty a_n z^{n-1},\,\,\, |z| < 1$$
    where $(a_n)$ is a sequence of complex numbers satisfying the following two conditions:
     $$\sum_{n=1}^\infty |a_n|^2 < \infty \mbox{ and } \int_0^1 |f(e^{2\pi i t})|^2 \varphi(t) dt < \infty.$$
In fact it is shown in \cite{Mukeru_Mulaudzi} that the space $\mathscr{H}$ is a Hilbert space with respect to the inner product
      $$\langle f, g \rangle = \int_0^1 f(e^{2\pi i t}) \overline{g(e^{2\pi i t})} \varphi(t) dt.$$
      (Here $f(e^{2\pi i t})$ is the left limit $\lim_{r \uparrow 1} f(r e^{2\pi i t})$.)
Moreover $\mathscr{H}$ is a reproducing kernel Hilbert space associated to the kernel
 $$\mathbb{K}: D \times D \to \mathbb{C},  \,\,\,   \mathbb{K}(z, w) = \sum_{n,m=1}^\infty \left(G^{-1}\right)_{k,j} z^{k-1} (\overline{w})^{j-1}.$$
 This means that for every function $f\in \mathscr{H}$ and for all $|z| < 1$, 
     $$ f(z) = \int_0^1 f(e^{2\pi i t})\overline{\mathbb{K}(e^{2\pi i t}, z)} \varphi(t)dt.$$ Or equivalently,
   $$ f(z) = \int_0^1 f(e^{2\pi i t})\mathbb{K}(z, e^{2\pi i t}) \varphi(t)dt$$ because obviously 
   $$\overline{\mathbb{K}(y, z)}  = \mathbb{K}(y, z)$$ by the fact that $G$ is hermitian. We shall write
   $$\mathbb{K}(z, w) = \sum_{n,m=1}^\infty \left(G^{-1}\right)_{k,j} z^{k-1} (\overline{w})^{j-1} = Z^T G^{-1}\overline{W}$$
   where $$Z^T = (1, z, z^2, z^3, \ldots)\mbox{ and } W^T = (1, w, w^2, w^3, \ldots).$$
   
A key property of a reproducing kernel Hilbert space is that the quantity $\sum_{n=1}^\infty \phi(z) \overline{\phi(w)}$ is invariant for all possible choices of orthonormal basis $\phi_1, \phi_2,\phi_3,\ldots$  and in fact it is equal to the kernel of the space. 
Now it is again proven in \cite{Mukeru_Mulaudzi} that the monomials $f(z) = z^{n-1}$  for all $n = 1,2,\ldots$ are elements of $\mathscr{H}$ and indeed the sequence $1, z, z^2, z^3, \ldots$ forms a basis in $\mathscr{H}$ (but in general it is not an orthonormal basis). This implies that any infinite sequence of polynomials $Q_1(z), Q_2(z), Q_3(z), \ldots$ such that $\mbox{deg}Q_n(z) = n-1$ for all $n = 1,2,3\ldots$ is a basis of $\mathscr{H}$. In particular, the sequence $P_1(z), P_2(z), \ldots$ (of the previous section) is an orthonormal basis of $\mathscr{H}$. Since this particular orthonormal basis is such that
$$\sum_{n= 1}^\infty P_n(z) \overline{P_n(w)} =\frac{ \psi(z) \overline{\psi(w)}}{1 - z \overline{w}}$$ it follows 
that the kernel $\mathbb{K}$ of $\mathscr{H}$ is such that:
 $$\mathbb{K}(z, w) = \sum_{n= 1}^\infty P_n(z) \overline{P_n(w)} =\frac{ \psi(z) \overline{\psi(w)}}{1 - z \overline{w}},\,\,\,\mbox{ for all } |z|<1, |w| < 1.$$
That is,
\begin{eqnarray}
Z^T G^{-1} \overline{W} = \frac{ \psi(z) \overline{\psi(w)}}{1 - z \overline{w}}.
\end{eqnarray}
 This is a key result that will now be used to compute the matrix $G^{-1}$. 
 
 \section{Explicit expression of the inverse matrix $G^{-1}$}
Since 
 \begin{eqnarray} \label{SM1}
Z^T G^{-1} \overline{W} = \frac{ \psi(z) \overline{\psi(w)}}{1 - z \overline{w}}
\end{eqnarray}
for all $z,w$ in the unit disc, then  taking $z = w = 0$ yields that $Z^T = W^T = (1,0,0,0, \ldots)$ and hence
$$Z^T G^{-1} \overline{W} = (G^{-1})_{1,1}$$ that is the element of the first row and first column of $G^{-1}$. Hence
      $$ (G^{-1})_{1,1} =  \psi(0) \overline{\psi(0)} = |\psi(0)|^2 = e^{- \int_0^1 \log(\varphi(t)) dt}. $$
To obtain $\left(G^{-1}\right)_{2,1}$ one can take differentiate (\ref{SM1}) with respect to $z$ at the point $z = w  = 0$ and obtain 
  $$\frac{\partial}{\partial z} \left(Z^T G^{-1} \overline{W}\right)|_{z=w=0} = (G^{-1})_{2,1}.$$ Hence
$$(G^{-1})_{2,1} = \left. \frac{\partial}{\partial z} \left(\frac{ \psi(z) \overline{\psi(w)}}{1 - z \overline{w}}\right)\right|_{z=w=0} = \psi'(0) \overline{\psi(0)}.$$ 
Also 
         $$(G^{-1})_{3,1} = {\scriptstyle\frac{1}{2}}\frac{\partial^2}{\partial z^2} \left(Z^T G^{-1} \overline{W}\right)|_{z=w=0} = {\scriptstyle\frac{1}{2}}\left. \frac{\partial^2}{\partial z^2} \left(\frac{ \psi(z) \overline{\psi(w)}}{1 - z \overline{w}}\right)\right|_{z=w=0} = {\scriptstyle\frac{1}{2}} \psi''(0) \overline{\psi(0)}.$$
 In general for all $k = 1,2,3, \ldots$
   $$(G^{-1})_{k,1} = \frac{1}{(k-1)!} \psi^{(k-1)}(0) \overline{\psi(0)}.$$
 Since the matrix $G^{-1}$ is also hermitian, then clearly
  $$(G^{-1})_{1,j} = \frac{1}{(j-1)!}\overline{\psi^{(j-1)}(0)} \psi(0).$$
 All other elements can also be obtained by taking successive partial derivatives with respect to $z$ and $\overline{w}$, that is, 
\begin{eqnarray} \label{main_SM_formula}
(G^{-1})_{k,j} = \left. \frac{1}{(k-1)! (j-1)!}\frac{\partial^{k+j-2}}{\partial z^{k-1}\partial (\overline{w})^{j-1}} \left(\frac{ \psi(z) \overline{\psi(w)}}{1 - z \overline{w}}\right)\right|_{z=w=0}.
\end{eqnarray}
For example,
 $$(G^{-1})_{2,2} = \left. \frac{\partial^{2}}{\partial z \partial \overline{w}} \left(\frac{ \psi(z) \overline{\psi(w)}}{1 - z \overline{w}}\right)\right|_{z=w=0}  = |\psi(0)|^2 + |\psi'(0)|^2.$$
 Also 
\begin{eqnarray*}
(G^{-1})_{3,3} & = & |\psi(0)|^2 + |\psi'(0)|^2 + \frac{1}{4} |\psi''(0)|^2,\\
 (G^{-1})_{4,4} & = & |\psi(0)|^2 + |\psi'(0)|^2 + \frac{1}{4} |\psi''(0)|^2 + \frac{1}{36} |\psi^{(3)}(0)|^2,\\
(G^{-1})_{5,5} & = & |\psi(0)|^2 + |\psi'(0)|^2 + \frac{1}{4} |\psi''(0)|^2 + \frac{1}{36} |\psi^{(3)}(0)|^2 + \frac{1}{576} |\psi^{(4)}(0)|^2.
\end{eqnarray*}  
An induction argument yields that
$$(G^{-1})_{n,n} = \sum_{k=0}^{n-1} \frac{|\psi^{(k)}(0)|^2}{(k!)^2}$$ and hence
 \begin{eqnarray} \label{Muk02}
\lim_{n\to \infty} (G^{-1})_{n,n} = \sum_{k=0}^{\infty} \frac{|\psi^{(k)}(0)|^2}{(k!)^2}.
\end{eqnarray}
Since  the Szeg\"o function $S(z)$ is analytic in the unit disc $D =\{z \in \mathbb{C}: |z| < 1\}$ and the Szeg\"o condition implies that it does not vanish anywhere in $D$, then its inverse $\psi(z) = 1/S(z)$ is also an analytic function in $D$. Assume then that
 $$\psi(z) = \sum_{n=0}^\infty a_n z^n,\,\,z\in D,\,\, a_n \in \mathbb{C}.$$
Then 
$$\frac{\psi^{(k)}(0)}{(k!)} = a_k$$ and hence (\ref{Muk02}) yields 
 \begin{eqnarray} \label{Muk03}
\lim_{n\to \infty} (G^{-1})_{n,n} =  \sum_{k=0}^\infty |a_k|^2.
\end{eqnarray}
In particular $\sum_{k=0}^\infty |a_k|^2 < \infty$. By the Szeg\"o theorem, since $1/\varphi(t)$ is also integrable and satisfies Szeg\"o condition, then 
      $$\int_0^1 \frac{e^{-2\pi i n t} }{\varphi(t)} dt = \sum_{k=0}^\infty \overline{a_n} a_{n+k}$$
for all $k = 0,1,2,\ldots$
 Taking $k = 0$ yields
            $$\int_0^1 \frac{1}{\varphi(t)} dt = \sum_{k=0}^\infty |a_n|^2,$$
  and therefore      
 $$\lim_{n\to \infty} (G^{-1})_{n,n} =  \sum_{k=0}^\infty |a_k|^2 = \int_0^1 \frac{1}{\varphi(t)} dt = \widehat{1/\varphi}(0).$$
 In general taking $$\psi(z) = \sum_{n=0}^\infty a_n z^n$$ in the main formula (\ref{main_SM_formula}) yields that
 for all $k, j$ with $j \leq k$,
   \begin{eqnarray} \label{sdsdeddsw344}
(G^{-1})_{k,j} = \overline{a_0} a_{\ell} + \overline{a_1} a_{\ell+1} + \overline{a_2} a_{\ell+2} + \ldots + \overline{a_{j-1}} a_{k-1},\,\,\, \mbox{where } \ell = k-j,\, j \leq k.
\end{eqnarray}
Thus the $G^{-1}$ is fully determined by the coefficients of the function $\psi(z)$. Hence 
 \begin{eqnarray*} 
\lim_{n \to \infty} (G^{-1}){k+n, k} =  \sum_{\ell=0}^\infty \overline{a_\ell} a_{n+\ell} = \int_0^1 \frac{e^{-2\pi n i t}}{\varphi(t)} dt = \widehat{1/\varphi}(-n).
\end{eqnarray*}
Hence we retrieve the following well-known result (Whittle approximation):
\begin{eqnarray*} 
\lim_{k \to \infty} (G^{-1}){k+n, k} = \int_0^1 \frac{e^{-2\pi n i t}}{\varphi(t)} dt = \widehat{1/\varphi}(-n).
\end{eqnarray*}


\section{The upper left-hand bloc of $G^{-1}$}
From the discussion above, it is clear that once the elements of the first row of $G^{-1}$ are determined, the other elements can easily be determined. 
Let $\left(G^{-1}\right)_{n\times n}$ be the upper left-hand $n \times n$ bloc of $G^{-1}$ (that is the bloc of $G^{-1}$ consisting of the first $n$ rows and first $n$ columns of $G^{-1}$). Then relation (\ref{sdsdeddsw344}) implies that all the elements of $\left(G^{-1}\right)_{n\times n}$ are explicitly determined by the first $n-1$ derivatives of the function $\psi(z)$ at the origin (or the first $n$ coefficients of $\psi(z)$ in its Taylor expansion at zero). In fact, 
\begin{eqnarray} \label{dsdsw334rer}
\left(G^{-1}\right)_{n\times n} = 
\begin{pmatrix} \overline{a_0} & 0 & 0 & 0 & \ldots & 0\\
                \overline{a_1} & \overline{a_0} & 0 & 0  &\ldots & 0\\
                 \overline{a_2} & \overline{a_1} & \overline{a_0} & 0 & \ldots & 0\\
                 \vdots & \vdots & \vdots & \vdots &  & \vdots\\
                 \overline{a_{n-1}} & \overline{a_{n-2}} & \overline{a_{n-3}} & \overline{a_{n-4}} & \ldots & \overline{a_0}             
 \end{pmatrix}
 \begin{pmatrix} a_0 & a_1 & a_2 & a_3& \ldots & a_{n-1}\\
                0 & a_0 & a_1 & a_2  &\ldots & a_{n-2}\\
                 0 & 0 & a_0 & a_1 & \ldots & a_{n-3}\\
                 \vdots & \vdots & \vdots & \vdots &  & \vdots\\
                 0 & 0 & 0 & 0 & \ldots & a_0             
 \end{pmatrix}
\end{eqnarray}
This is yields an LU decomposition of the matrix $G^{-1}$. This decomposition also yields a system of orthonormal polynomials with respect to the weight function $1/\varphi(t)$. Indeed it is immediately seen that 
the system 
 $\{Q_0(z), Q_1(z), \ldots, Q_{n}(z), \ldots\}$ given by:
  $$Q_n(z) =  a_n  + a_{n-1} z + a_{n-2} z^2 + \ldots + a_0 z^n,\,\,\,\, z\in D, n = 0,1,2,\ldots$$
 is orthonormal on the unit circle with respect to the weight function $1/\varphi(t)$ (the inverse of the spectral density of $G$). 
 
\section{Illustrating example: fractional Gaussian noise}
For the fractional Gaussian noise of index $0 < H < 1$, as already discussed, 
  \begin{eqnarray*} 
\psi_H(z) = \exp\left({-\scriptstyle\frac{1}{2}}\int_0^1 \left(\frac{e^{2\pi i t} + z}{e^{2\pi i t} - z}\right) \log(\varphi_H(t)) dt\right),\,\, \mbox{ for all } |z| < 1  
\end{eqnarray*}
  and 
 \begin{eqnarray*}
\varphi_H(t)= (\cos(2\pi t)-1)\left(\mbox{Li}(-2 H, e^{-2\pi i t}) + \mbox{Li}(-2 H, e^{2\pi i t})\right),\,\,0 < t < 1, 
\end{eqnarray*}
where $\mbox{Li}$ is the analytic continuation of the classical polylogarithm function 
 $$\mbox{Li}(s,z) = \sum_{k=1}^\infty  k^{-s} z^k.$$ 
We can therefore calculate the derivatives of $\psi_H(z)$ at the origin with respect to $z$. Note that 
 $$\frac{d}{dz} \left(\int_0^1 \left(\frac{e^{2\pi i t} + z}{e^{2\pi i t} - z}\right) \log(\varphi_H(t)) dt\right) = \int_0^1 \frac{d}{dz} \left(\frac{e^{2\pi i t} + z}{e^{2\pi i t} - z}\right) \log(\varphi_H(t)) dt. $$ 
Set $$u_0 = {-\scriptstyle\frac{1}{2}}\int_0^1  \log(\varphi_H(t)) dt$$ and 
$$u_k = -\int_0^1 e^{-2\pi i k t} \log(\varphi_H(t)) dt,\,\, k = 1,2,\ldots.$$
 Then we can express the quantities $a_n = \psi^{(n)}(0)/(n!)$ as functions of $u_0, u_1, u_2, \ldots, u_n$. 
We shall write
   $$\psi(z) = e^{U(z)},\,\, U(z) = {-\scriptstyle\frac{1}{2}}\int_0^1 \left(\frac{e^{2\pi i t} + z}{e^{2\pi i t} - z}\right) \log(\varphi_H(t)) dt$$
   Then 
   $$\psi'(z) = U'(z) \psi(z).$$ 
   Using Leibniz formula this yields
$$\psi^{(n+1)}(0) = \sum_{k=0}^{n} \binom{n}{k}U^{(k+1)}(0) \psi^{(n-k)}(0)$$
where $$U^{(0)}(0) = U(0) = {-\scriptstyle\frac{1}{2}}\int_0^1  \log(\varphi_H(t)) dt$$ and for $k = 1,2,\ldots,$,
$$U^{(k)}(0) = - k!\int_0^1 e^{-2\pi i k t} \log(\varphi_H(t)) dt.$$
Then 
 $$ U^{(k)}(0) = - k!\int_0^1 e^{-2\pi i k t} \log(\varphi_H(t)) dt = k! u_{k},\,\,k = 1,2,\ldots$$
Hence
       $$a_{n+1} = \frac{\psi^{(n+1)}(0)}{(n+1)!} = \frac{1}{n+1}\sum_{k=0}^{n} (k+1) u_{k+1} a_{n-k}$$
with $$a_0 = \psi(0) = e^{u_0}.$$ 
These relations  yield the values of the coefficients $(a_n)$ for all $n$ and finally we can compute the bloc $\left(G^{-1}\right)_{n\times n}$ for any number $n$. 
Taking for example $n = 5$ and the Hurst index $H = 0.75$ yields that the 5 values of sequence $(u_k)$ are: 
\begin{eqnarray*}
u_0 =  0.113994, u_1 = -0.333504, u_2 = -0.123701, u_3 = -0.0838558, u_4 = -0.0626411
\end{eqnarray*}
 and the corresponding first 5 values of $(a_k)$ are:
 
\begin{eqnarray*}
 a_0 = 1.12075, a_1 = -0.373773, a_2 = -0.0763097, a_3 = -0.0546738, a_4 = -0.0374192.
\end{eqnarray*}
Hence this yields using relation (\ref{dsdsw334rer}) that the upper left-hand $5 \times 5$ bloc of $G^{-1}$  (that is the bloc corresponding of the first 5 rows and 5 columns) is: 
\begin{eqnarray*}
\left(G^{-1}\right)_{5\times 5} = 
\begin{pmatrix}
1.25607 &  -0.418904 &  -0.0855238 &  -0.0612754 & -0.0419375\\
-0.418904 &  1.39578 & -0.390382 & -0.0650882 & -0.0472891\\
-0.0855238 & -0.390382 & 1.4016 & -0.386209 & -0.0622327\\
-0.0612754 & -0.0650882 &  -0.386209 &  1.40459 & -0.384164\\
-0.0419375 & -0.0472891 & -0.0622327 & -0.384164 & 1.40599
\end{pmatrix}.
\end{eqnarray*} 
One can compare that with the inverse of the bloc $(G)_{m\times m}$ of the original covariance matrix $G$ for $m$ sufficiently large. Taking $m = 1000$ yields the following: 
\begin{eqnarray*}
\left(\left[\left(G\right)_{m\times m}\right]^{-1}\right)_{5 \times 5} = 
\begin{pmatrix}
1.25599 &  -0.418956 &  -0.0855708 &  -0.0613187 &  -0.0419781\\
-0.418956 & 1.39574 &  -0.390413 & -0.0651171 & -0.0473163\\
-0.0855708 & -0.390413 &  1.40157 &  -0.386235 &  -0.0622571\\
-0.0613187 &  -0.0651171 &  -0.386235 &  1.40456 & -0.384186\\
-0.0419781 &  -0.0473163 & -0.0622571 & -0.384186 & 1.40597
\end{pmatrix}.
\end{eqnarray*} 
These two matrices are very close and the equality holds when $m \to \infty$.

\section{Illustrating example: finite-diagonal Toeplitz matrices}
In this section we now consider the particular case where the Toeplitx matrix $G$ has only a finite number $2m+1$ of non-zero diagonals (such matrices are also called $2m+1$--diagonal Toeplitz matrices). It is also an interesting problem to compute the inverse of such matrices. 
Assume as previously that 
  $G = (g_{k,j})$ with 
  $$g_{k,j} = \gamma(j-k) = \overline{\gamma(k-j)},\,\,k,j = 1,2,3,\ldots$$ and $\gamma(0) = 1$. Moreover we assume that $\gamma(k) = 0 $ for the integers $|k|\geq m+1.$This means that the matrix $G$ is a $2m+1$--diagonal Toeplitz matrix. Set $$\gamma(k) = q_k \mbox{ for all } k = 1,2,3,\ldots, m.$$
The first row of $G$ is:
  $$1, q_1, q_2, \dots, q_m, 0, 0, 0, \ldots$$
The corresponding spectral density if obviously given by
 $$\varphi(t) = 1 + \sum_{k=1}^m e^{2k\pi i t} q_k +  \sum_{k=1}^me^{-2 k\pi i t} \,\overline{q_k},\,\,\, t\in [0,1].$$
The condition that $G$ is positive definite implies that $\varphi(t) \geq 0$  everywhere everywhere in $[0, 1]$. If the function $\log(\varphi(t))$ is integrable, then the corresponding Szeg\"o function $S(z)$ and its inverse $\psi(z)$ are given as already discussed by:
  \begin{eqnarray*} 
S(z) = \exp\left({\scriptstyle\frac{1}{2}}\int_0^1 \left(\frac{e^{2\pi i t} + z}{e^{2\pi i t} - z}\right) \log(\varphi(t)) dt\right),\,\, \mbox{ for all } |z| < 1.  
\end{eqnarray*}
  \begin{eqnarray*} 
\psi(z) = S(z)^{-1} = \exp\left({-\scriptstyle\frac{1}{2}}\int_0^1 \left(\frac{e^{2\pi i t} + z}{e^{2\pi i t} - z}\right) \log(\varphi(t)) dt\right),\,\, \mbox{ for all } |z| < 1.  
\end{eqnarray*}
The formulas of the previous section can be used to compute the elements $(a_n)$ and from which  the inverse matrix $G^{-1}$ follows.  
\subsection{Case of infinite 3--diagonal Toeplitz matrix}
This corresponds to $m = 1$. Assume that the matrix  $G$ is defined by the function
\begin{eqnarray*}
\gamma(k) = \left\{\begin{array}{cc}
         1 & \mbox{ if } k = 0\\
         q & \mbox{ if } k = 1\\
         \overline{q} & \mbox{ if } k = -1\\
         0 &  \mbox{ otherwise}
         \end{array}\right.
         \end{eqnarray*}
where $q$ is a complex number such that $|q| < 1/2.$ The corresponding spectral density is 
 $$\varphi(t) = 1 + q \,e^{2\pi i t} + \overline{q}\, e^{-2\pi i t},\,\,t\in [0, 1].$$     
We can explicitly compute the corresponding function $\psi(z)$ as follows. First consider the Szeg\"o function $S(z)$ associated to $\phi(t)$ and write as in (\ref{szego1}):
Write 
$$S(z) = c_0 + c_1 z + c_2 z^2 + c_3 z^3 + \ldots$$ where $c_0, c_1, c_2,\ldots$ are complex numbers. We can compute the values of these parameters of $S(z)$ using similar argument as for the parameters $a_0, a_1, a_2, \ldots$. 
Because 
      $$S(z) = \psi(z)^{-1} = \exp\left({\scriptstyle\frac{1}{2}}\int_0^1 \left(\frac{e^{2\pi i t} + z}{e^{2\pi i t} - z}\right) \log(\varphi(t)) dt\right)$$
Then 
\begin{eqnarray} \label{esdwdsada12}
c_{n+1} = \frac{S^{(n+1)}(0)}{(n+1)!} = \frac{1}{n+1}\sum_{k=0}^{n} (k+1) v_{k+1} c_{n-k}\end{eqnarray}
with 
$c_0 = S(0) $ and        
$$v_k = - u_k = \int_0^1 e^{-2\pi i k t} \log(\varphi(t)) dt,\,\, k = 1,2,\ldots.$$      
This yields that    
\begin{eqnarray*}
c_0 = \sqrt{\frac{1  + \sqrt{1 - 4|q|^2}}{2}} \mbox{ and } c_1 = \frac{q}{c_0}
\end{eqnarray*}
and $c_k = 0$ for all $k = 2,3,\ldots$
Hence the Szeg\"o function associated to $\varphi$ is:
  $$S(z) = c_0 + \frac{q}{c_0} z,\,\, |z| < 1.$$
  The corresponding function $\psi$ is:
\begin{eqnarray} \label{SMKKA1}
\psi(z) =\frac{1}{S(z)}  = \frac{c_0}{c_0^2 + q z}\mbox{ where } c_0 = \sqrt{\frac{1  + \sqrt{1 - 4|q|^2}}{2}}.
\end{eqnarray}
That is:
       $$\psi(z) = \frac{\sqrt{(1  + \sqrt{1 - 4|q|^2})/2}}{((1  + \sqrt{1 - 4|q|^2})/2) + q z}.$$
In the case where $q$ is a real number, this function was obtained in \cite{Mukeru_Mulaudzi} where it is expressed as:
        $$\psi(z) = \left(\frac{2}{|q|}\right)^{1/2} \left(\frac{1}{a + b z}\right)\mbox{ where } a = \sqrt{|q|^{-1} + \sqrt{|q|^{-2} - 4}} \mbox{ and } b = (2/a) \mbox{ sign}(q).$$
We can now explicitly compute the inverse of the matrix $G$.    
Using (\ref{SMKKA1}), it is clear that the Taylor expansion of $\psi(z)$ is given by: 
     $$\psi(z) = \sum_{n=0}^\infty (-1)^n \left(\frac{q^n}{c_0^{2n+1}}\right) z^n$$ and hence the coefficients $a_0, a_1, a_2, \ldots$ of $\psi(z)$ are explicitly given by:
     $$a_n = (-1)^n \left(\frac{q^n}{c_0^{2n+1}}\right) ,\,\,n = 0,1,2, \ldots.$$
Therefore the element of order $(k, j)$ of the inverse matrix $G^{-1}$ (with $j \leq k$) is explicitly given by:
 $$(G^{-1})_{k,j}  = \sum_{\ell=1}^{j} \overline{a_{k-\ell}}\,\, a_{j-\ell}  = \frac{(-1)^{j+k} c_0^{2(1-k-j)} (\overline{q})^{k-j} \left(c_0^{4 j} - |q|^{2j}\right)}{c_0^4 - |q|^2}$$ where 
   $$c_0 = \sqrt{\frac{1  + \sqrt{1 - 4|q|^2}}{2}}.$$
 This yields an explicit analytic expression of the inverse of the matrix $G$. 
Summarising we have the following result:
 \begin{proposition}
For any complex number $q$ such that $|q| < 1/2$, the infinite matrix $G$ given by
 \begin{eqnarray*}
G_{k,j} = \left\{\begin{array}{cc}
         1 & \mbox{ if } k = j\\
         q & \mbox{ if } j- k = 1\\
         \overline{q} & \mbox{ if } j-k = -1\\
         0 &  \mbox{ otherwise}
         \end{array}\right.
         \end{eqnarray*}
is such that its inverse is the  matrix $G^{-1}$ given explicitly by:
\begin{eqnarray*}
\left(G^{-1}\right)_{k,j} =  \frac{(-1)^{j+k} c_0^{2(1-k-j)} (\overline{q})^{k-j} \left(c_0^{4 j} - |q|^{2j}\right)}{c_0^4 - |q|^2},\,\,\,\, \mbox{ for }j \leq k 
\end{eqnarray*} 
and$$G_{k,j} =  \overline{\left(G^{-1}\right)_{j,k}}\,\,\,  \mbox{ for  } j > k $$
   where 
   $$c_0 = \sqrt{\frac{1  + \sqrt{1 - 4|q|^2}}{2}}.$$    
\end{proposition}
To illustrate this result, take $q = -1/5$. Then the upper left-hand $5 \times 5$ bloc of $G^{-1}$ is:
\begin{eqnarray*}
\left(G^{-1}\right)_{5\times 5} = 
\begin{pmatrix}
\frac{5 (5 - \sqrt{21})}{2} & \frac{ 25(5 - \sqrt{21})}{2} & \frac{550- 120\sqrt{21}}{2}  & \frac{5(527 - 115 \sqrt{21})}{2} &  \frac{5(2525 - 551\sqrt{21})}{2}  \\

\frac{ 25(5 - \sqrt{21})}{2} & \frac{ 25(23 - \sqrt{21})}{2} & \frac{2750 - 600 \sqrt{21}}{2}  & \frac{25(527 - 115\sqrt{21})}{2} & \frac{25(2525 - 551\sqrt{21})}{2} \\

 \frac{550- 120\sqrt{21}}{2} & \frac{2750 - 600 \sqrt{21}}{2}& \frac{ 13200 - 2880\sqrt{21}}{2} & \frac{63240 - 13800\sqrt{21}}{2} & \frac{120(2525 - 551  \sqrt{21})}{2}\\
 
 \frac{5(527 - 115 \sqrt{21})}{2} & \frac{25(527 - 115\sqrt{21})}{2}   & \frac{63240 - 13800\sqrt{21}}{2} & \frac{575(527 - 115 \sqrt{21})}{2} &  \frac{575(2525 - 551 \sqrt{21})}{2}  \\
 
  \frac{5(2525 - 551\sqrt{21})}{2}  &  \frac{25(2525 - 551\sqrt{21})}{2}  & \frac{120(2525 - 551  \sqrt{21})}{2} & \frac{575(2525 - 551 \sqrt{21})}{2} & \frac{2755 (2525 - 551 \sqrt{21})}{2}
\end{pmatrix}
\end{eqnarray*}
We can again compare this exact matrix with the inverse of the bloc $(G)_{m\times m}$ of  $G$ for $m$ sufficiently large. Taking $m = 10$ yields the following: 
\begin{eqnarray*}
\left(\left[\left(G\right)_{m\times m}\right]^{-1}\right)_{5 \times 5} = 
\begin{pmatrix}
1.04356 &  0.217804 &  0.0454583 &  0.0094877 & 0.0019802\\
0.217804 & 1.08902 &  0.227292 & 0.0474385 & 0.00990099\\
0.0454583 &  0.227292 & 1.091 &  0.227705 &  0.0475248\\
0.0094877 &  0.0474385 &  0.227705 & 1.09109 &  0.227723\\
0.0019802 &  0.00990099 &  0.0475248 &  0.227723 &  1.09109
\end{pmatrix},
\end{eqnarray*}
which is very close to the exact matrix. 
\subsection{Case of infinite 5--diagonal Toeplitz matrices}
Assume $G$ is given by the function:
 \begin{eqnarray*}
\gamma(k) = \left\{\begin{array}{cc}
         1 & \mbox{ if } k = 0\\
         q_1 & \mbox{ if } k = 1\\
          q_2 &\mbox{ if } k = 2\\
         \overline{q_1} & \mbox{ if } k = -1\\
           \overline{q_2} & k = -2 \\
         0 &  \mbox{ otherwise}
         \end{array}\right.
         \end{eqnarray*}
Write 
       $$\varphi(t) = 1 +  e^{2\pi i t}\, q_1 + e^{-2\pi i t}\, \overline{q_1} +  e^{4\pi i t}\, q_2 + e^{-4\pi i t}\, \overline{q_2}.$$
Write
  $$S(z) = c_0 + c_1 z + c_2 z^2 + \ldots$$ and we obtain that 
    $$c_0 = \exp\left({\scriptstyle\frac{1}{2}}\int_0^1 \log(\varphi(t)) dt\right)$$
    $$c_1 = \frac{c_0(c_0^2 q_1 - q_2 \overline{q_1})}{c_0^4 - |q_2|^2},\,\,\,\, c_2 = q_2/c_0 \mbox{ and } c_3 = c_4 = \ldots = 0.$$
In particular case where $q_1, q_2$ are real numbers, then 
 $$c_1 = \frac{c_0 q_1}{c_0^2 + q_2}.$$
Therefore the Szeg\"o function $S(z)$ is the polynomial:
 $$S(z) = c_0 + c_1 z + c_2 z^2 =  c_0 + \left(\frac{c_0(c_0^2 q_1 - q_2 \overline{q_1})}{c_0^4 - |q_2|^2}\right) z + \left(\frac{q}{c_0}\right) z^2,\,\, |z| < 1.$$         
So we only need to compute one integral to obtain the value of $c_0$. We do not know a closed form for the integral giving $c_0$ so one will need to a numerical approximation of $c_0$ which is very easy to obtain. 
Finally the inverse Szeg\"o function $\psi(z)$ is:
 $$\psi(z) = \frac{1}{S(z)} = \frac{c_0(c_0^4 - |q_2|^2)}{c_0^2(c_0^4 - |q_2|^2) + c_0^2(c_0^2 q_1 - q_2 \overline{q_1}) z + q(c_0^4 - |q_2|^2) z^2}.$$
We can now easily express $\psi(z)$ in the form $\psi(z) = a_0 + a_1 z + a_2 z^2 + a_3 z^3 + \ldots$ and deduce the elements of the inverse matrix $G^{-1}$.\\ 
To illustrate consider the following infinite matrix $G = (\gamma(j-k))_{k,j}$ with 
 \begin{eqnarray*}
\gamma(k) = \left\{\begin{array}{cc}
         1 &\mbox{ if } k = 0\\
         -1/4 & \mbox{ if } |k| = 1\\
          1/3 &\mbox{ if } |k| = 2\\
           0 &  \mbox{ otherwise}
         \end{array}\right.
         \end{eqnarray*}
The spectral density function corresponding to $G$ is:
 $$\varphi(t) = 1 -  (e^{2\pi i t} +e^{-2\pi i t})/4+ (e^{4\pi i t} + e^{-4\pi i t})/3. $$
 It is easy to check that $\varphi(t)$ is positive everywhere and hence $G$ is positive definite. 
 Then 
 $$ c_0 = 0.909567,\,\,  c_1 =  -0.195918, c_2 =  0.366475 \mbox{ and }  c_3 = c_4 = \ldots = 0.$$
Hence
      $$S(z) = 0.909567  -0.195918 z + 0.366475 z^2$$ and 
      $$\psi(z) = S(z)^{-1} = (0.909567  -0.195918 z + 0.366475 z^2)^{-1}.$$ 
For example the first 5 coefficients in the Taylor expansion of $\psi(z)$ at $z = 0$ are"
 $$a_0 = 1.09942, a_1 = 0.236813, a_2 = -0.391961, a_3 = -0.179842, a_4 = 0.119188.$$ This yields that the the upper left-hand $5\times 5$ bloc of $G^{-1}$ is: 
\begin{eqnarray*}
\left(G^{-1}\right)_{5\times 5} = 
\begin{pmatrix} 
 1.20873& 0.260358& -0.430932 & -0.197723& 0.131038\\
0.260358 &   1.26481 & 0.167536 & -0.473521 & -0.169497\\
-0.430932 & 0.167536 & 1.41845 & 0.238028 & -0.520238\\
-0.197723 & -0.473521 & 0.238028 & 1.45079 & 0.216593\\
0.131038 & -0.169497 & -0.520238 &  0.216593 & 1.465
\end{pmatrix}.
\end{eqnarray*}
We can again compare with the inverse of the bloc $(G)_{m\times m}$ of  $G$ for $m$ sufficiently large. Taking $m = 100$ yields the following: 
\begin{eqnarray*}
\left(\left[\left(G\right)_{m\times m}\right]^{-1}\right)_{5 \times 5} = 
\begin{pmatrix}
1.20873& 0.260358& -0.430932 & -0.197723& 0.131038\\
0.260358 &   1.26481 & 0.167536 & -0.473521 & -0.169497\\
-0.430932 & 0.167536 & 1.41845 & 0.238028 & -0.520238\\
-0.197723 & -0.473521 & 0.238028 & 1.45079 & 0.216593\\
0.131038 & -0.169497 & -0.520238 &  0.216593 & 1.465
\end{pmatrix}.
\end{eqnarray*}
(This is exactly the same as the previous matrix and hence, the error of approximation is $\leq 10^{-6}$.)

\subsection{General case }
In the  general case of $2m+1$--diagonal Toeplitz matrix, write
      $$\varphi(t) = 1 + \sum_{k=1}^m e^{2k\pi i t} q_k +  \sum_{k=1}^m e^{-2\pi i t} \,\overline{q_k}$$ where $ q_1, q_2,\ldots, q_m$ are complex numbers. In line with the results of the previous two particular cases ($m = 1$ and $m = 2$), one can conjecture that the corresponding Szeg\"o function is a polynomial of degree $\leq m$, i.e.
 $$S(z) = c_0 + c_1 z + c_2 z^2 + \ldots + c_m z^m,\,\,z \in D$$ where $c_0, c_1, \ldots, c_m$ are  complex numbers. Then under that assumption one only needs  to compute these coefficients from relations (\ref{esdwdsada12}) and deduce that
 $$\psi(z) = \frac{1}{S(z)} = \frac{1}{c_0 + c_1 z + c_2 z^2 + \ldots + c_m z^m}.$$
For example, consider  the matrix 
  $G = (\gamma(j-k))_{k,j}$ with 
   \begin{eqnarray*}
\gamma(k) = \left\{\begin{array}{cc}
         1 & \mbox{ if } k = 0\\
         3/10 & \mbox{ if } k = 1\\
          2 (1+ i)/10 &\mbox{ if } k = 2\\
          (1 + i)/10 &\mbox{ if } k = 3\\
         3/10 & \mbox{ if } k = -1\\
         2(1-i)/10 &\mbox{ if } k = -2\\
           (1-i)/10 &\mbox{ if } k = -3\\
         0 &  \mbox{ otherwise}
         \end{array}\right.
         \end{eqnarray*}
  
Then 
 \begin{eqnarray*}
\varphi(t) & = & (1/10)\left(10 + 3 e^{2\pi i t} + (2 + 2i)e^{4\pi i t} + (1 + i) e^{6\pi i t} + 3 e^{-2\pi i t} + (2 - 2i)e^{-4\pi i t} + (1 - i) e^{-6\pi i t}\right)\\
& = & (1/5) (5 + 3 \cos(2\pi t) + 2 \cos(4 \pi t) + \cos(6 \pi t) - 
   2 \sin(4 \pi t) - \sin(6 \pi t))
\end{eqnarray*}
We compute
   $$c_0 =e^{{\scriptstyle\frac{1}{2}}\int_0^1 \log(\varphi(t)) dt} = 0.917429$$
and 
 $$c_1 = 0.242589 - 0.0634194 i ,\,\, c_2 = 0.196713 + 0.181643i,\,\, c_3 = 0.109 + 0.109 i$$ 
and $c_n = 0$ for all $n \geq 4$.  
Then 
    $$S(z) = 0.917429 + (0.242589 - 0.0634194 i) z + (0.196713 + 0.181643i)z^2 + (0.109 + 0.109 i)z^3$$ and 
  $$\psi(z) = (0.917429 + (0.242589 - 0.0634194 i) z + (0.196713 + 0.181643i)z^2 + (0.109 + 0.109 i)z^3)^{-1}.$$
From this explicit inverse Szeg\"o function $\psi(z)$ one can now compute every element of the inverse matrix $G^{-1}$. For example using (\ref{main_SM_formula}), we have that 
     $$\left(G^{-1}\right)_{10,9} = -0.282433 +  0.183806\, i.$$
 We can also compute the first 10 coefficients $a_0, a_1, \ldots, a_9$ of the Taylor expansion of $\psi(z)$ at $z = 0$ and  using (\ref{sdsdeddsw344}) obtain
$$\left(G^{-1}\right)_{10,9} = \sum_{k=0}^8 \overline{a_{k}}\,\, a_{k+1} = -0.282433 + 0.183806\,i.$$
The first $5\times 5$ bloc of $G^{-1}$ is: 
 $$ \left(G^{-1}\right)_{5\times 5} = A + i B$$ where 
\begin{eqnarray*}
 A =  
\begin{pmatrix} 
1.18811 & -0.314162  & -0.177357  &  
  0.00862465  &  0.0300867\\
-0.314162 & 1.27685&  -0.286529  &  -0.182067 &
  0.00981616\\
-0.177357  &  -0.286529  & 1.36869&  -0.279574  &  -0.217594 \\
0.00862465  &  -0.182067 & -0.279574  &  
  1.36979 &  -0.283269 \\
0.0300867  & 0.00981616  &  -0.217594  & -0.283269  &  1.38529
\end{pmatrix}
\end{eqnarray*}
and \begin{eqnarray*}
 B =  
\begin{pmatrix} 
0&  0.0821306  &  - 0.278669 &  
  - 0.0351419  &    0.132323 \\
 - 0.0821306  & 0 &   0.168077  &  - 0.269973 &
 - 0.0722108 \\
 0.278669  &  - 0.168077  & 0 &    0.175346  &  - 0.282669  \\
 0.0351419  &   0.269973  & - 0.175346  &  
 0 & 0.177196 \\
- 0.132323  &  0.0722108  &   0.282669  &  - 
   0.177196  &  0
\end{pmatrix}.
\end{eqnarray*}
 Again this exact value is very close to the first $5 \times 5$ bloc of the inverse of $\left(G\right)_{m\times m}$ for $m$ large. In fact $m=100$ yields an approximation error $\leq 10^{-6}$. 

\section{Concluding remark} 
This paper studies the inverse of  infinite Toeplitz matrices $G$ that are hermitian and positive definite. We have obtained an exact formula that expresses the elements of the inverse in terms of the coefficients of the inverse Szeg\"o function associated to the matrix $G$. We provided explicit calculations throughout in order to  clearly show how the method can be used in practice to compute the inverse of an infinite matrix. We applied the results to the covariance matrix of the classical fractional Gaussian noise and  this is the first time  the inverse of this matrix is explicitly computed. This will be of interest in the study of the inverse problem of fractional Brownian motion. We  also considered the particular case of infinite Toeplitz matrices with only a finite number $2m+1$ of non-zero diagonals. For the case $m = 1$ (tridiagonal matrices) we obtained an explicit analytical formula for the inverse. We also considered the case where $m = 2$ (pentadiagonal matrices) but here no explicit analytic solution is known. An interesting observation for the inverse of $2m+1$--diagonal Toeplitz matrices is that,  possibly, the associated Szeg\"o function is a polynomial of the form $S(z) = c_0 + c_1 z + c_2 z^2 + \ldots + c_m z^m$.  This means the following: 
  If $q_1, q_2, \ldots, q_m$ are complex numbers such that the function
    $$\varphi(t) = 1 + \sum_{k=1}^m e^{2k\pi i t} q_k +  \sum_{k=1}^me^{-2\pi i t} \,\overline{q_k}$$ is strictly positive for all $t \in [0, 1]$, then there exist complex coefficients $c_1, c_2, \ldots, c_m$ such that 
  $$ \exp\left({\scriptstyle\frac{1}{2}}\int_0^1 \left(\frac{e^{2\pi i t} + z}{e^{2\pi i t} - z}\right) \log(\varphi(t)) dt\right)  = c_0 + c_1 z + c_2 z^2 + \ldots + c_m z^m$$
  for all complex numbers $z$ such that $|z| < 1$.  
This needs to be proven. In that case the inverse of the matrix is obtained by simply determining the Taylor coefficients $(a_0, a_1, a_2, \ldots)$ of $(c_0 + c_1 z + c_2 z^2 + \ldots + c_m z^m)^{-1}$ and making use of formula (\ref{dsdsw334rer}).


\begin{thebibliography}{100}



\bibitem{Beran} Beran, J.  1994. {\it  Statistics for Long-Memory Processes.} Chapman and Hall.

\bibitem{Caflisch} Caflisch, R.E. 1981. 
An inverse problem for Toeplitz matrices and the synthesis of discrete transmission lines. {\it Linear Algebra and its Applications,} 38, 207--225

 \bibitem{Cintoli_et_al} Cintoli, S., Shlomo P. Neuman, S.P. and  Di Federico, V. 2005.  Generating and scaling fractional Brownian motion on finite domains, {\it Geophysical Research Letters}, Vol.32, L08404.
 


 
 
 \bibitem{Dambrogi-ola} D'Ambrogi-Ola, B. 2009. {\it Inverse problem of fractional Brownian motion with discrete data}, PhD Thesis, University of Helsinki.
 


\bibitem{Delorme_2015} Delorme, M.  and Wiese, K.J. 2015. The maximum of a fractional
Brownian motion: analytic results from perturbation theory. {\it Phys. Rev. Lett. } 115,  210601. 


 

 \bibitem{daFonseca} da Fonseca, C.M. and  Petronilho, J. 2001.  Explicit inverses of some tridiagonal matrices
{\it Linear Algebra and its Applications} 325,  7--21.



 


 \bibitem{Szego} Grenander, U. and Szeg\"o, G. 1958. {\it Toeplitz forms and their applications.} University
of Calif. Press, Berkeley and Los Angeles. 

 

\bibitem{Mukeru_Pisieral} Mukeru, S. 2020. A generalisation of Pisier homogeneous Banach algebra. To appear in {\it Michigan Math. J.}

\bibitem{Mukeru_Mulaudzi} Mukeru, S. and Mulaudzi, M.P. 2020. Zeros of Gaussian power series, Hardy spaces and  determinantal point processes, {\it Preprint}

\bibitem{Pedro} Mentz, R.P. On the inverse of some covariance matrices of Toeplitz type. {\it SIAM J. Appl. Math.} 31(3), 426--437
 
\bibitem{Rostek} Rostek, S. and Sch\"obel, R. 2013. A note on the use of fractional Brownian motion for financial modeling. {\it Economic Modelling} 30, 30--35.


\bibitem{Simon1} Simon, B. 2009.  Orthogonal Polynomials on the Unit Circle, AMS Colloquium Series, vol. 1, American Mathematical Society, Rhode Island


  \bibitem{Sinai} Sinai,Y.G. 1976. Self-similar probability distributions. {\it Theory of Probability and its Applications} 21, 64--80. 
  
 \bibitem{Vehel} V\'ehel, J.L. and Riedi, R. 1997. Fractional Brownian motion and data traffic modeling: the other end of  the spectrum. {\it Fractals in Engineering,} 97, 185--202.
 
\bibitem{Wang} Wang, C., Li, H., and Zhao, D. 2015. An explicit formula for the inverse of a pentadiagonal Toeplitz matrix. {\it J. Comput. Appl. Math.}  278, 12--18.

\bibitem{Whittle} Whittle, P. 1953.  Estimation and information in stationary time series. {it Ark.
Mat.} 2, 423-434.

\end{thebibliography}
\end{document}